\documentclass[12pt]{article}

\usepackage{amssymb}
\usepackage{amsmath}
\usepackage{amsbsy}
\usepackage{amscd}
\usepackage{amsfonts}
\usepackage{amsthm}
\usepackage{mathrsfs}
\usepackage{verbatim}
\usepackage[colorlinks]{hyperref}
\usepackage{fullpage}
\usepackage{mathdots}
\usepackage{graphicx,subfigure}
\usepackage[english]{babel}
\usepackage[utf8]{inputenc}
\usepackage{tikz}
\usepackage{adjustbox}
\usepackage{authblk}
\usepackage{caption}

\usepackage{mathtext}

\usepackage{tikz-cd}
\usetikzlibrary{backgrounds,fit, matrix}
\usetikzlibrary{positioning}
\usetikzlibrary{calc,through,chains}
\usetikzlibrary{arrows,shapes,snakes,automata, petri}
\usepackage[boxsize=1.25em, centerboxes]{ytableau}

\newtheorem{Theorem}[equation]{Theorem}

\newtheorem{Lemma}[equation]{Lemma}
\newtheorem{Proposition}[equation]{Proposition}

\theoremstyle{definition}

\newtheorem{Example}[equation]{Example}

\newtheorem{Notation}[equation]{Notation}

\newtheorem{Remark}[equation]{Remark}

\numberwithin{equation}{section}
\numberwithin{figure}{section}

\newcommand{\F}{{\mathbb F}}

\newcommand{\Z}{{\mathbb Z}}
\newcommand{\Q}{{\mathbb Q}}
\newcommand{\R}{{\mathbb R}}

\newcommand{\mb}[1]{\mathbf{#1}}
\newcommand{\mc}[1]{\mathcal{#1}}

\newcommand{\vep}{\varepsilon}

\DeclareMathOperator{\pyr}{pyr}

\DeclareMathOperator{\conv}{conv}

\DeclareMathOperator{\HH}{\mathbf{H}}

\begin{document}

\title{Toric Codes from Order Polytopes}

\author[1]{Mahir Bilen Can}
\author[2]{Takayuki Hibi}

\affil[1]{\small{
Tulane University, New Orleans, Louisiana, US\\
mahirbilencan@gmail.com}}

\affil[2]{\small{Department of Pure and Applied Mathematics, Graduate School of Information Science and Technology, Osaka University, Suita, Osaka 565--0871, Japan
hibi@math.sci.osaka-u.ac.jp}}

\maketitle

\begin{abstract}
In this article we investigate a class of linear error correcting codes in relation with the order polytopes.
In particular we consider the order polytopes of tree posets and bipartite posets. 
We calculate the parameters of the associated toric variety codes. 
\vspace{.5cm}


\end{abstract}

\section{Introduction}

In the present article we are concerned with a special class of algebraic-geometric codes~\cite{TVN} that are defined on toric varieties. 
Building on a work of S. Hansen~\cite{Hansen2001}, J. Hansen initiated the study of toric codes on polygons in~\cite{Hansen2002}.
This development quickly led to numerous new results on the algebraic-geometric codes that are constructed on higher dimensional toric varieties. 
The articles~\cite{LS2006,LS2007, Ruano2007,SS2010} amplified the importance of combinatorial approach in determining the parameters of the toric codes. 
Our goal in this article is to show that, the set of order polytopes form an interesting ground for the applications of such work.

Let $P$ be a poset whose elements are listed as $\vep_1,\dots, \vep_m$. 
Let $N$ denote the free $\Z$-module on $P$, $N:=\bigoplus_{i=1}^m \Z \varepsilon_i$.
Let $M$ denote the dual of $N$, that is $M:=Hom_\Z (N,\Z)$.
The dual of the element $\vep_i$ ($i\in \{1,\dots, m\}$) in $M$ will be denoted by $e_i$. 
Let $2^P$ denote the set of all subsets of $P$. We define the function $\rho : 2^P \to N\otimes_\Z \Q$ by $W \mapsto \sum_{\varepsilon_i \in W} \vep_i$.
The {\em order polytope} of $P$, denoted by $\mb{O}_P$, is the convex hull of the finite set 
\[
\{ \rho(W):\ \text{ $W$ is an upper order ideal of $P$}\}.
\] 
The face lattice of the polytope $\mb{O}_P$ was first described by Geissinger~\cite{Geissinger1981}, 
whose results were amplified by Stanley in~\cite{Stanley1986}.
A concrete description of the edges of $\mb{O}_P$ can be found in~\cite{HLSS}.
Following~\cite{Hibi1987}, we now introduce a class of toric varieties that are closely related to the order polytopes. 
The set of all order ideals of $P$, denoted by $J(P)$, is a distributive lattice with respect to inclusion.
In particular, we have the joins (denoted by $\vee$) and the meets (denoted by $\wedge$) of the elements of $J(P)$. 
Let $Y:=\{y_\alpha :\ \alpha \in J(P)\}$ be a set of algebraically independent variables indexed by the order ideals.
Then the {\em Hibi toric scheme associated with $P$} is the projective scheme $\text{Proj} \ k[Y]/I$, where $I$ is the homogeneous ideal 
\[
I = ( y_\alpha  y_\beta - y_{\alpha \wedge \beta} y_{ \alpha \vee \beta} :\ y_\alpha, y_\beta \in Y).
\]
It turns out that the fan of $X_P$ is the normal fan of the order polytope $\mb{O}_P$.

The purpose of our article is to investigate the parameters of the toric code of the defining polytope $\mb{O}_P$ of $X_P$.
The parameters that we speak of are called the ``length,'' the ``dimension,'' and the ``minimum distance.'' 
Although our method applies to all finite posets, in this article we focus on the minimum distance computation for the order polytopes of the rooted trees only. Let $P=\{\vep_1,\dots, \vep_m\}$ be a rooted tree, where $\vep_1$ is the root. 
We view $P$ as a connected, graded poset with the unique minimal element as the root.
Our first main result (recorded as Theorem~\ref{T:mainresult}) states that minimum distance of the toric code $\mc{C}_{\mb{O}_P}$ over a finite field $\F_q$, where $q> 3$, is given by 
\[
d(\mc{C}_{\mb{O}_P}) = (q-1)^a (q-2)^b,
\]
for some $a$ and $b$ such that $a+b=m$. In fact, we know precisely what $a$ and $b$ are. 

Let $\mb{P}$ be a polytope. 
The {\em length} of the associated toric code $\mc{C}_{\mb{P}}$ over $\F_q$ is given by $(q-1)^{\dim \mb{P}}$, 
where ${\dim \mb{P}}$ is the dimension of the affine hull of $\mb{P}$. 
Hence, in our case, the length is given by $(q-1)^{\dim \mb{O}_P} = (q-1)^m$, where $m$ is the cardinality of the poset $P$.  
On the other hand, the {\em dimension} of a toric code of $\mb{P}$ is given by the number of lattice points in $\mb{P}$. 
Therefore, in our case, it is given by the number of (upper) order ideals of $P$. 
For a rooted tree with $m$ vertices, this number (dimension) varies in the range $m+1,\dots, 2^{m-1}+1$; 
it is equal to the number of order preserving maps $\sigma : P \to \{ 0,1\}$. 
The unique rooted tree with $m$ vertices that has $m+1$ order ideals is the chain with $m$ vertices. 
The unique rooted tree with $m$ vertices that has $2^{m-1}+1$ order ideals is the ``$m$-th shrub'' defined in Section~\ref{S:Main1}.

Let $Q$ be a graded poset with $2m$ elements ($m\in \Z^+$).
If $Q$ has $m$ minimum elements, then we will call $Q$ an {\em $(m,m)$-bipartite poset}. 
The second infinite family of toric codes that we consider comes from the order polytopes of $(m,m)$-bipartite posets.
Our second main result (recorded as Theorem~\ref{T:mainresult2}) states that the minimum distance of the toric code $\mc{C}_{\mb{Q}}$ over a finite field 
$\F_q$ where $q>3$ is given by 
\[
d(\mc{C}_{\mb{O}_P}) = (q-1)^m (q-2)^m. 
\]
The dimension of such a code varies in the range $2^{m+1}-1,\dots, 3^m$.

Before closing this introduction, we want to mention a fact we inferred from our calculations. 
In general, a preferable linear error correcting code is the one that has a ratio of $\text{dimension}/\text{length}$ fixed while the ratio 
$\text{minimum distance}/\text{length}$ is as large as possible.  
It is natural to wonder if it is possible to increase these ratios for a toric code by switching to the polar polytope. 
In this article we pay a close attention to the polar of the order polytope of a graded poset. 
It turns out that, by a result of Hibi and Higashitani~\cite{HH2011}, the polar polytope of a suitable dilation of $\mb{O}_P$, called the {\em poset polytope of $P$}, is reflexive and terminal. 
(We will explain these notions in the sequel.) 
These properties essentially imply that the number of lattice points of a poset polytope is much smaller compared to the number of lattice points of the order polytope. 
Hence, as far as the parameters of linear codes are concerned, the order polytopes are better than the poset polytopes.

The structure of our paper is as follows.
In the next section we introduce our basic notation regarding posets, polytopes, and toric codes. 
In the same section we briefly review  some results of Soprunov and Soprunova also. 
The purpose of Section~\ref{S:HH} is to compare the structures of the order polytopes and poset polytopes. 
We prove our first main result about the toric codes defined by the rooted tree posets in Section~\ref{S:Main1}.
We prove our second main result about the toric codes defined by the $(m,m)$-bipartite graphs in Section~\ref{S:Main2}.
In addition, in this section, we observe that (Lemma~\ref{L:ordinalfree}) the free sum of two order polytopes, $\mb{O}_P\oplus \mb{O}_Q$, 
is equivalent to the order polytope $\mb{O}_{P\oplus Q}$, where $P\oplus Q$ stands for the ordinal sum of $P$ and $Q$. 
Here, the equivalence is defined by the change of coordinates.

\section{Preliminaries} \label{S:Preliminaries}

In this article, by a {\em poset} we will always mean a finite poset.  
A {\em lower order ideal} in $P$ is a subposet $I$ such that for every $y\in I$,
if $x\leq y$ in $P$, then $x\in I$. 
An {\em upper order ideal} in $P$ is defined similarly where we replace the condition $x\leq y$ with $y\leq x$.

The set of all lower order ideals of $P$ is denoted by $J(P)$. 
This is a distributive lattice with respect to inclusion. 
The set of all upper order ideals also form of a distributive lattice, which is isomorphic to $J(P^{opp})$,
where $P^{opp}$ denotes the opposite poset to $P$. 
An order reversing bijection between two posets will be called an {\em anti-isomorphism}. 
If $P$ and $Q$ are two isomorphic (resp. anti-isomorphic) posets, then we will write $P\cong Q$ (resp. $P\cong_a Q$).

Let $x$ and $y$ be two elements from $P$. 
If $x \leq y$, and $x\leq z \leq y$ implies that $z=x$ or $z = y$, then $y$ is said to {\em cover} $x$.  
Customarily, the cover relation is denoted by $x\lessdot y$.

A {\em chain} is a poset $C:=\{x_1,\dots, x_n\}$ whose elements are linearly ordered, $x_1 \lneq x_2 \lneq \cdots \lneq x_n$.
A {\em maximal chain} in a poset $P$ is a chain $C\subseteq P$ such that $C$ is not a subposet of any other chain in $P$.
If $C=\{x_1,\dots, x_k\}$ is a chain, then the {\em length of $C$} is defined as $k-1$.

An {\em antichain} is a poset whose elements are all incomparable. 
The greatest possible size of an antichain in $P$ is called the {\em width of $P$}. 
Dilworth's theorem~\cite{Dilworth} states that the width is equal to the minimal number of chains that cover the set.

A poset $P$ is called a {\em graded (or ranked) poset} if every maximal chain in $P$ has the same length. 
In this case, a function $\ell: P \to \Z$ which has the property that 
$\ell(y) = \ell(x) +1$ for every cover relation $x\lessdot y$ in $P$ is called a {\em rank function} for $P$.
Without loss of generality we assume that $\ell(x)=0$ whenever $x$ is a minimal element. 
Then $\ell$ is uniquely determined by $P$, so, we call it {\em the} rank function of $P$.

The {\em Hasse diagram} of a poset $P$ is the directed graph whose vertex set is the set of elements of $P$
such that for $x,y\in P$ there is a directed edge from $x$ to $y$ if $x$ is covered by $y$ in $P$. 
A poset $P$ is said to be {\em connected} if its Hasse diagram is connected. 
Clearly, if a finite poset possesses a top element (denoted by $\hat{1}$) or a bottom element (denoted by $\hat{0}$), then it is connected. 
A {\em lattice} is a poset $L$ such that every pair of elements has a least upper bound and a greatest lower bound. 

The {\em polar} (or {\em dual}) of a polytope $\mb{P}\subset \Q^m$ is the polytope $\mb{P}^\circ$ defined by
\[
\mb{P}^\circ:= \{ y\in (\Q^m)^*:\ \text{$\langle x,y \rangle \leq 1$ for all $x \in \mb{P}$}\}.
\]
Here, $\langle, \rangle$ is the canonical evaluation pairing between $\Q^m$ and $(\Q^m)^*$.

Let $x_0$ be a point in $\Q^{m}$, and let $H$ be a hyperplane in $\Q^m$ such that $x_0\notin H$. 
Let $\mb{P}$ be a polytope in $H$. 
The {\em pyramid over $\mb{P}$ with apex at $x_0$} is the convex hull $\text{conv}(\mb{P}, x_0)$.
We will denote a pyramid over $\mb{P}$ by $\pyr(\mb{P})$.

The vertex set of a polytope $\mb{P}$ will be denoted by $V(\mb{P})$.
Let $\mb{Q}$ and $\mb{P}$ be two polytopes in $\Q^m$ and $\Q^n$, respectively. 
The {\em direct product} (or simply the {\em product}) of $\mb{Q}$ and $\mb{P}$, denoted by $\mb{Q}\times \mb{P}$, is defined as the convex hull, 
\begin{align*}
\mb{Q}\times \mb{P}:= \text{conv} ( (a,b):\ a\in V(\mb{Q}), \ b\in V(\mb{P})).
\end{align*}
We now assume that the origin of $\Q^m$ (resp. of $\Q^n$) is contained in $\mb{Q}$ (resp. in $\mb{P}$).
The {\em free sum of $\mb{Q}$ and $\mb{P}$}, denoted by $\mb{Q}\oplus \mb{P}$, is defined as follows: 
\begin{align*}
\mb{Q}\oplus \mb{P}:= \text{conv} ( \mb{Q}\times \{0_{\Q^n}\}, \{0_{\Q^m}\} \times \mb{P}).
\end{align*}

\subsection{Toric codes.}

The purpose of this subsection is to introduce {\em toric codes} by circumventing much of the original definition of the algebraic-geometric codes. For a detailed introduction to this important subject, we recommend the textbook~\cite{TVN}. 

Let $N$ be a free abelian group of rank $m$, and let $M$ denote its dual group. 
Let $\mb{P}$ be a full dimensional lattice polytope in $M\otimes_\Z \Q$.
The lattice points in $\mb{P}\cap M$ define monomials that are regarded as polynomial functions on the $m$-dimensional torus $T_N:= \text{Hom}(N,\overline{\F_q^*})$. 
Let $H^0(T_N(\F_q),\mb{P})$ denote the $\F_q$-vector space that is spanned by these monomials. 
The {\em toric code of $\mb{P}$} is then the image of the evaluation map 
\begin{align*}
\text{ev}: H^0(T_N(\F_q), \mb{P}) &\longrightarrow (\F_q^*)^m \\
f &\longmapsto (f(x))_{x\in T_N(\F_q)}. 
\end{align*}
More generally, the algebraic-geometric code associated with an ample line bundle on a normal variety $X$ that is defined over $\overline{\F_q}$
is the image of the germ-evaluation map on a set of $\F_q$-rational points $S\subseteq X(\F_q)$. 
The toric codes from lattice polytopes are defined by evaluating on the $\F_q$-rational points of the open orbit of a normal toric variety. 

Hereafter, we denote by $\mc{C}_{\mb{P}}$ the toric code associated with a lattice polytope $\mb{P}$. 
The {\em length of $\mc{C}_{\mb{P}}$} is defined as 
\[
\text{length } := (q-1)^m,
\]
where $m$ is the dimension of the toric variety. 
The {\em dimension of $\mc{C}_{\mb{P}}$} is defined as the vector space dimension of the space of sections 
\[
\text{ dimension } := \dim H^0(T_N(\F_q), \mb{P}).
\]
This number is given by the number of lattice points $\mb{P}\cap M$.
Finally, the computation of the minimum distance for the toric codes associated with an order polytope is the main focus of the present article. 
It is calculated as follows. 
For a section $f \in H^0(T_N(\F_q),\mb{P})$, let $Z(f)$ denote the number of points in $(\F_q^*)^m$ where $f$ vanishes. 
Then the {\em minimum distance of $\mc{C}_{\mb{P}}$}, denoted by $d(\mc{C}_{\mb{P}})$, is given by 
\begin{align*}
d(\mc{C}_{\mb{P}})= (q-1)^m - \max_{f\in H^0(T_N(\F_q),\mb{P})\setminus \{0\}} Z(f).
\end{align*}

We will make use of the following results which are due to Soprunov and Soprunova. 
\begin{Lemma}(Theorem 2.1~\cite{SS2010})\label{L:SS2010T2.1}
Let $\mb{P}$ and $\mb{Q}$ be two lattice polytopes contained in the boxes $[0,q-2]^m\subseteq \Q^m$ and $[0,q-2]^n\subseteq \Q^n$, respectively. 
Then the minimum distance of the code of the product $\mb{P}\times \mb{Q}$ is given by 
$d(\mc{C}_{\mb{P}\times \mb{Q}}) = d(\mc{C}_{\mb{P}}) d(\mc{C}_{\mb{Q}})$.
\end{Lemma} 
 
Let $K_q^n$ denote the $n$-dimensional cube $[0,q-2]^n$. Let $\mb{Q}$ be an $n$-dimensional lattice polytope contained in $K_q^n$. Then the {\em unit pyramid over $\mb{Q}$} is defined by $\text{conv} \{ e_{n+1}, (x,0): x\in \mb{Q} \}$, where $e_{n+1}$ is the unit vector $(0,\dots, 0,1)\in \R^{n+1}$.

\begin{Lemma}(Theorem 2.3~\cite{SS2010})\label{L:SS2010T2.3}
Let $\mb{Q}$ be a lattice polytope of $\dim \mb{Q}\geq 1$.
If $\mb{P}$ denotes the unit pyramid over $\mb{Q}$, then we have  
$d(\mc{C}_{\mb{P}}) = (q-1)  d(\mc{C}_{\mb{Q}})$.
\end{Lemma}

\section{Order Polytopes, Poset polytopes}\label{S:HH}

Let $P= \{ \vep_1,\dots, \vep_m \}$ be a finite poset, and let $N$ denote the free $\Z$-module generated by $P$.
Let $\hat{P}$ denote $P \cup \{ \hat{0},\hat{1}\}$, where $\hat{0}$ (resp. $\hat{1}$) is such that $\hat{0} \lneq \vep_i$ (resp.
$\vep_i \lneq \hat{1}$) for every $i\in \{1,\dots, m\}$. 
Let $M$ denote the dual of $N$, that is $M:=Hom_\Z (N,\Z)$, and let $\{e_1,\dots, e_m\}$ be the basis of $M$ that is dual to $P$. 
Let us temporarily denote $\hat{0}$ (resp. $\hat{1}$) by $\vep_0$ (resp. $\vep_{m+1}$).
Then for each covering relation $\vep_i \lessdot \vep_j$ in $\hat{P}$, we introduce a vector $\rho(\vep_i,\vep_j)$ in $M\otimes_\Z \Q$ as follows:
\begin{align}\label{A:linearform}
\rho(\vep_i,\vep_j):= 
\begin{cases}
e_i & \text{ if $\vep_j= \hat{1}$;}\\
e_i-e_j & \text{ if $\vep_i,\vep_j \in P$;}\\
-e_j & \text{ if $\vep_i=\hat{0}$.}
\end{cases}
\end{align}
The {\em poset polytope of $P$}, denoted by $\HH_P$, is the convex hull of points $\rho(\vep_i,\vep_j)$, where $\vep_i\lessdot \vep_j$ is a cover in $\hat{P}$.
A systematic study of these polytopes is initiated by Hibi and Higashitani in~\cite{HH2011}.
In this article, we construct linear error correcting codes by using (the polars of the) poset polytopes.

Next, we will discuss poset polytopes and their relationship to the order polytopes.
Since it is already introduced (in the Introduction), we will not repeat the definition of a poset polytope here. 
In~\cite{HH2011}, Hibi and Higashitani showed that these polytopes have some remarkable properties. 
We will summarize the relevant results from~\cite{HH2011} in the form of a single lemma to ease our referencing.

\begin{Lemma} \label{L:HH}
For every poset $P$, the following statements hold:
\begin{enumerate}
\item $\HH_P$ is a {\em Fano polytope}, that is, 0 is the unique integral interior point. 
\item $\HH_P$ is {\em terminal}, that is, each integral point on the boundary of $\HH_P$ is a vertex.
\item $\HH_P$ is {\em Gorenstein}, that is, its dual polytope is integral.
\item If $P$ is a graded poset of length $l-2$, then the polar polytope of $\HH_P$ is the dilated and translated order polytope $l \mb{O}_P-v$, where $v$ is the unique lattice point in $l\mb{O}_P$.
\end{enumerate}
\end{Lemma}

The first item is proved in~\cite[Lemma 1.3]{HH2011}, 
the second item is proved in~\cite[Lemma 1.4]{HH2011}. 
The third item is proved in~\cite[Lemma 1.5]{HH2011}. 
The last item is recorded in~\cite[Remark 1.6]{HH2011}; its proof follows from the definitions. 

\begin{Remark}
A Gorenstein and Fano polytope is known as the {\em reflexive polytope}. 
In particular, the dual of a reflexive polytope is reflexive. 
The normal fan of a reflexive polytope gives a ``Gorenstein Fano toric variety''~\cite[Theorem 8.3.4]{CLS}. 
(Such toric varieties are always normal.)
In particular, a reflexive polytope is very ample in the sense of~\cite[Definition 2.2.17]{CLS}.
\end{Remark}

\begin{Notation}
If $P$ is a graded poset of length $l-2$, then the polytope $l \mb{O}_P-v$, where $v$ is the unique lattice point in $l\mb{O}_P$, will be denoted by $\mb{O}_P(l)$.
\end{Notation}

\begin{Example}\label{E:vposet1}

Let $P$ (resp. $\hat{P}$) be the poset whose Hasse diagram is on the left (resp. on the right) in Figure~\ref{F:2}.
\begin{figure}[htp]
\begin{center}
\scalebox{.9}{
\begin{tikzpicture}[scale=1]
\begin{scope}[xshift=-3.5cm]
\node at (0,-1.75) {$P$};
\node (1) at (0,0) {$\varepsilon_1$};
\node (2) at (-1,1) {$\varepsilon_2$};
\node (3) at (1,1) {$\varepsilon_3$};
\draw[thick,-] (1) to (2);
\draw[thick,-] (1) to (3);

\end{scope}
\begin{scope}[xshift=3.5cm]
\node at (0,-1.75) {$\hat{P}$};
\node (0) at (0,-1) {$\hat{0}$};
\node (1) at (0,0) {$\varepsilon_1$};
\node (2) at (-1,1) {$\varepsilon_2$};
\node (3) at (1,1) {$\varepsilon_3$};
\node (4) at (0,2) {$\hat{1}$};
\draw[thick,-] (0) to (1);
\draw[thick,-] (1) to (2);
\draw[thick,-] (1) to (3);
\draw[thick,-] (2) to (4);
\draw[thick,-] (3) to (4);
\end{scope}
\end{tikzpicture}
}
\end{center}
\caption{}
\label{F:2}
\end{figure}

By fixing $\{\vep_1,\vep_2,\vep_3\}$ as a basis for $N\otimes_\Z \Q$, we will identify the elements of $N\otimes_\Z \Q$ by their coordinate vectors.
Then, the vertex set of $\mb{O}_P$ consists of the following vectors in $\Q^3$:
\begin{align*}
\rho(\emptyset) &= (0,0,0), \\
\rho(\{\varepsilon_2\}) &= \vep_2 = (0,1,0), \\
\rho(\{\varepsilon_3\}) &= \vep_3 = (0,0,1), \\
\rho(\{\varepsilon_2,\varepsilon_3\}) &= \vep_2+ \vep_3 = (0,1,1), \\
\rho(\{\varepsilon_1,\varepsilon_2,\varepsilon_3\}) &= \vep_1+ \vep_2+ \vep_3 =(1,1,1).
\end{align*}
In Figure~\ref{F:thepolytope}, we depicted the order polytope of $P$. 
\begin{figure}[htp]
\begin{center}
\scalebox{.9}{
\begin{tikzpicture}[scale=2]
\filldraw[fill=green!20,draw=green!50!black, ultra thick] (0,0) -- (0,1.45)-- (-1,.5) -- (-1,-1)  -- (0,0)  ;
\filldraw[fill=green!20,draw=green!50!black, ultra thick] (0,1.45) -- (.75,0.5) -- (-1,.5) -- (0,1.45)  ;
\filldraw[fill=green!20,draw=green!50!black, ultra thick] (.75,0.5) -- (-1,.5) -- (-1,-1) -- (.75,.5) ;
\draw[draw=green!50!black, ultra thick] (0,0) -- (.75,.5);
\draw[draw=green!50!black, ultra thick] (0,0) -- (-1,-1);
\draw[draw=green!50!black, ultra thick] (0,0) -- (0,1.45);
\draw[dashed, < - > ] (0,-1) -- (0,2); 
\draw[dashed, < - > ] (-1,0) -- (2,0); 
\draw[dashed, < - > ] (-1.75,-1.75) -- (.5,.5); 
\node at (0,0) {$\bullet$}; 
\node at (0,1.45) {$\bullet$}; 
\node at (-1,-1) {$\bullet$}; 
\node at (-1,0.5) {$\bullet$}; 
\node at (.75,0.5) {$\bullet$}; 
\node at (2.25,0) {$\vep_1$}; 
\node at (-1.4,-1.7) {$\vep_2$}; 
\node at (0,2.2) {$\vep_3$}; 
\node at (-.35,.1) {$(0,0,0)$}; 
\node at (-.35,1.55) {$(0,0,1)$}; 
\node at (-1.35,-1.1) {$(0,1,0)$}; 
\node at (-1.35,.55) {$(0,1,1)$}; 
\node at (1.1,.55) {$(1,1,1)$}; 
\end{tikzpicture}
}
\end{center}
\caption{The order polytope of $P$.}
\label{F:thepolytope}
\end{figure}
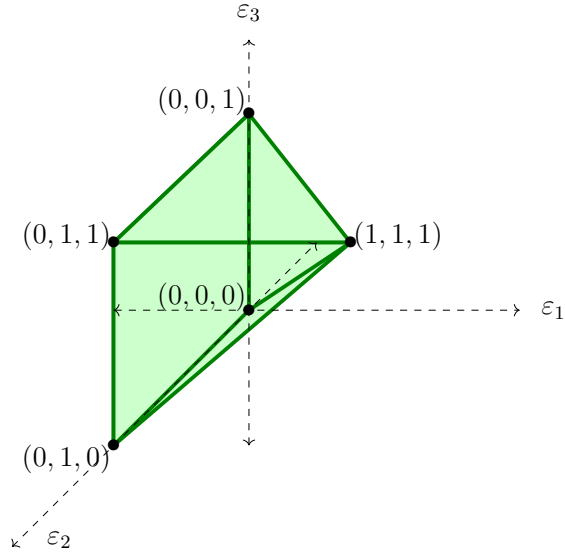
Finally, let us consider the dual polytope for $\mb{O}_P(3)$. 
It is easy to check that the vertices of the dual polytope $\HH_P$ are given by $-e_1,e_1-e_2,e_1-e_3, e_2,e_3$.
We notice that the convex hull of $e_1-e_2,e_1-e_3, e_2,e_3$ is a rectangular plate, which we denote by $A$. 
Then $\HH_P$ is a pyramid over $A$ with apex at $-e_1$.
\end{Example}

We close this subsection by two simple observations. 
\begin{Lemma}\label{L:freesum}
Let $P$ be a poset with connected components $P_1,\dots, P_r$. 
Then we have 
\[
\mb{H}_P = \mb{H}_{P_1} \oplus \cdots \oplus \mb{H}_{P_r}.
\]
\end{Lemma}
\begin{proof}
Let $x$ be a vertex in $\mb{H}_P$. 
Then there is a covering relation $\vep_i \lessdot \vep_j$ in $\hat{P}$ such that 
\begin{align*}
x\in \{ 
e_i,
e_i-e_j ,
-e_j \}.\end{align*}
Since every covering relation in $\hat{P}$ is a covering relation in one of the posets $\hat{P_i}$ ($i\in \{1,\dots, r\}$),
we see that the vertex set of $\mb{H}_P$ is a disjoint union,
\[
V(\mb{H}_P) = V(\mb{H}_{P_1})\sqcup \cdots \sqcup V(\mb{H}_{P_r}).
\]
Note that, the subpolytopes $\mb{H}_{P_i}$ for $i\in\{1,\dots, r\}$ are contained in skew subspaces in $\Q^m$.
Nevertheless, they all share the origin of $\Q^m$. 
Therefore, we have 
\begin{align*}
\mb{H}_{P} &= \conv(V( \mb{H}_P)) \\
&= \conv(V( \mb{H}_{P_1}) \sqcup \cdots \sqcup V(\mb{H}_{P_r})) \\ 
&= \conv(V( \mb{H}_{P_1})) \sqcup \cdots \sqcup  \conv(V(\mb{H}_{P_r})).
\end{align*}
This finishes the proof of our assertion. 
\end{proof}

Our next observation is about the order polytopes.
\begin{Lemma}\label{L:product}
Let $P$ be a poset with connected components $P_1,\dots, P_r$. 
Then we have 
\[
\mb{O}_P = \mb{O}_{P_1} \times \cdots \times \mb{O}_{P_r}.
\]
\end{Lemma}
\begin{proof}
Let $x$ be a vertex in $\mb{O}_P \subseteq \Q^m$, 
where $m$ is the number of elements of $P$.
Then there is an upper order ideal $I$ in $P$ such that $x=\rho(I)$. 
Since $P$ is the disjoint union $P_1\sqcup \cdots \sqcup P_r$, we see that $I= I_1 \sqcup \cdots \sqcup I_r$,
where $I_i$ ($i\in \{1,\dots, r\}$) is an upper order ideal in $P_i$. 
It follows that $x$ is of the form 
\begin{align}\label{A:Oproduct}
x= x_1+\dots+ x_r \in \Q^{m_1}\oplus \cdots \oplus \Q^{m_r},
\end{align}
where $x_i = \rho(I_i)$, and $\Q^{m_i}$ is the vector subspace of $\Q^m$ that is spanned by the basis vectors corresponding to the elements of $P_i$ ($i\in \{1,\dots, r\}$). 
The decomposition in (\ref{A:Oproduct}) shows that the vertex set of $\mb{O}_P$ is the product of the vertex sets of the order polytopes $\mb{O}_{P_i}$,
\[
V(\mb{O}_P) = V(\mb{O}_{P_1}) \times \cdots \times V(\mb{O}_{P_r}).
\]
This finishes the proof. 
\end{proof}

The decompositions that we observed in Lemmas~\ref{L:product} and~\ref{L:freesum} can be obtained from each other by induction 
and the well-known polarity correspondence between the free sums and direct products of polytopes.

\begin{Remark}
As we mentioned in the introduction, a desirable code is the one with a high {\em transmission rate}, that is, $\text{dimension}/\text{length}$. 
The construction of $\mb{H}_P$ uses the cover relations in $P$ whereas the construction of $\mb{O}_P$ uses all upper order ideals in $P$.
In general the vertices of the latter polytope are much more numerous.
Therefore, for a generic poset $P$, the transmission rate of $\mc{C}_{\mb{H}_P}$ is very small compared to the transmission rate of $\mc{C}_{\mb{O}_P}$.
\end{Remark}

\section{Shrubs} \label{S:Main1}

We begin with a reduction result.

\begin{Proposition}\label{P:dofproducts}
Let $P$ be a poset with $r$ connected components $P_1,\dots, P_r$. 
Let $q$ be a prime power such that $q> 2$. 
Then the minimum distance of the toric code $\mc{C}_{\mb{O}_P}$ is given by 
\[
d(\mc{C}_{\mb{O}_P})= d(\mc{C}_{\mb{O}_{P_1}}) \cdot \hdots \cdot d(\mc{C}_{\mb{O}_{P_r}}).
\] 
\end{Proposition}

\begin{proof}
We know from Lemma~\ref{L:product} that $\mb{O}_P$ decomposes as a direct product, 
\[
\mb{O}_P = \mb{O}_{P_1} \times \cdots \times \mb{O}_{P_r}. 
\]
By applying induction with Lemma~\ref{L:SS2010T2.1}, we see that $d(\mc{C}_{\mb{O}_P})= d(\mc{C}_{\mb{O}_{P_1}}) \cdot \hdots \cdot d(\mc{C}_{\mb{O}_{P_r}})$.
\end{proof}

Next, we focus on the connected posets.

\begin{Proposition}\label{P:pyramid}
Let $P=\{\vep_1,\dots, \vep_m\}$ be a connected poset with a unique minimal element, $\vep_1$. 
If $P'$ is the poset obtained from $P$ by removing $\vep_1$, then we have 
\begin{align*}
d(\mc{C}_{\mb{O}_P}) = (q-1) d(\mc{C}_{\mb{O}_{P'}}).
\end{align*}
\end{Proposition}

\begin{proof}
Since $\vep_1$ is the smallest element in $P$, the upper order ideal generated by $\vep_1$ is the whole poset $P$.  
In particular, all coordinates of the corresponding vertex $x_0:=\rho(P)$ in $\Q^m$ is 1,  
\[
x_0=(1,\dots, 1) \in \Q^m.
\]
For every other vertex $x=(a_1,\dots, a_m)$ of $\mb{O}_P$ such that $x\neq x_0$, we have $a_1=0$. 
This means that the line segment between vertices $x_0$ and $x$ is an edge of the polytope $\mb{O}_P$. 
(Note that this observation follows from~\cite[Lemma 1.1 (a)]{HLSS} as well.)  
It follows that $\mb{O}_P$ is a pyramid over $\mb{O}_{P'}$. 
Now, the rest of the proof follows from Lemma~\ref{L:SS2010T2.3}.
\end{proof}

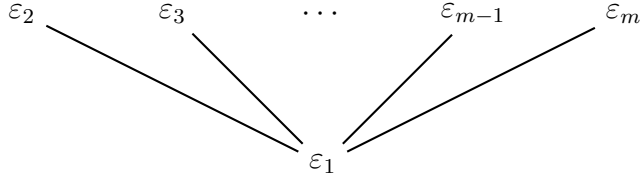
\begin{figure}[htp]
\begin{center}
\scalebox{1}{
\begin{tikzpicture}[scale=1]

\node (1) at (0,-1) {$\vep_1$};
\node (2) at (-4,1) {$\vep_2$};
\node (3) at (-2,1) {$\vep_3$};
\node (4) at (0,1) {$\cdots$};
\node (5) at (2,1) {$\vep_{m-1}$};
\node (6) at (4,1) {$\vep_m$};
\draw[thick,-] (1) to (2); 
\draw[thick,-] (1) to (3);
\draw[thick,-] (1) to (5); 
\draw[thick,-] (1) to (6); 
 
\end{tikzpicture}
}
\end{center}
\caption{The $m$-th shrub, $S_m$.}
\label{F:toric1}
\end{figure}

Let $P$ be a poset. 
We call $P$ a {\em rooted tree poset} if the following conditions hold: 
\begin{enumerate}
\item the Hasse diagram of $P$ is a rooted tree, where the smallest element of $P$ is the root;
\item the leaves of $P$ are the maximal elements of $P$.
\end{enumerate}
If $P$ is the rooted tree poset whose Hasse diagram is as in Figure~\ref{F:toric1}, then we call it the {\em $m$-th shrub}.
The $m$-th shrub will be denoted by $S_m$. 
If the number $m$ is understood from the context, or if it is not relevant to the discussion, then we simply write ``shrub'' instead of writing 
``the $m$-th shrub.'' 
Let $I$ be an upper order ideal in $S_m$. 
If $I$ contains the element $\vep_1$, then it is equal to $S_m$.
If $\vep_1 \notin I$, then $I$ can be any subset of $\{\vep_2,\dots, \vep_m\}$. 
Therefore, $J(S_m^{opp})$ is isomorphic to $B_{m-1}\oplus \hat{1}$, where $B_{m-1}$ is the boolean algebra of rank $m-1$. 
The proof of the following lemma is easy so we omit it.
\begin{Lemma}
Let $m\geq 2$. Then the order polytope of the shrub $S_m$ is a pyramid over the unit cube of dimension $m-1$. 
\end{Lemma}

\begin{figure}[htp]
\centering
\begin{tikzpicture}[scale=.75]

\begin{scope}[xshift=0cm]
\node at (0,0) {$\bullet$};
\node at (-4,1) {$\bullet$};
\draw[dotted, thick] (0,0) to  (-4,1);
\node at (-6,2) {$\bullet$};
\draw[dotted, thick] (-4,1) to  (-6,2);
\node at (-7,3) {$\bullet$};
\node at (-6,3) {$\bullet$};
\node at (-5,3) {$\bullet$};
\node at (-4,3) {$\bullet$};
\draw[-, thick] (-6,2) to  (-7,3);
\draw[-, thick] (-6,2) to  (-6,3);
\draw[-, thick] (-6,2) to  (-5,3);
\draw[-, thick] (-6,2) to  (-4,3);

\node at (-2,1) {$\bullet$};
\draw[dotted, thick] (0,0) to  (-2,1);
\node at (-2,2) {$\bullet$};
\node at (-2,3) {$\bullet$};
\draw[dotted, thick] (-2,1) to  (-2,2);
\draw[-, thick] (-2,2) to  (-2,3);
\node at (0,2) {$\bullet$};
\node at (-1,3) {$\bullet$};
\node at (0,3) {$\bullet$};
\node at (1,3) {$\bullet$};
\draw[dotted, thick] (0,0) to  (0,1);
\draw[dotted, thick] (0,1) to  (0,2);

\draw[-, thick] (0,2) to  (-1,3);
\draw[-, thick] (0,2) to  (0,3);
\draw[-, thick] (0,2) to  (1,3);

\node at (0,1) {$\bullet$};
\node at (2,3) {$\bullet$};
\node at (3,3) {$\bullet$};
\node at (4,3) {$\bullet$};
\node at (5,3) {$\bullet$};

\draw[dotted, thick] (0,0) to  (3,1);
\draw[dotted, thick] (3,1) to  (3,2);

\draw[-, thick] (3,2) to  (2,3);
\draw[-, thick] (3,2) to  (3,3);
\draw[-, thick] (3,2) to  (4,3);
\draw[-, thick] (3,2) to  (5,3);

\node at (3,1) {$\bullet$};
\node at (3,2) {$\bullet$};

\end{scope}
\end{tikzpicture}
\caption{The shrubbery of a tree.}
\label{F:Example}
\end{figure}
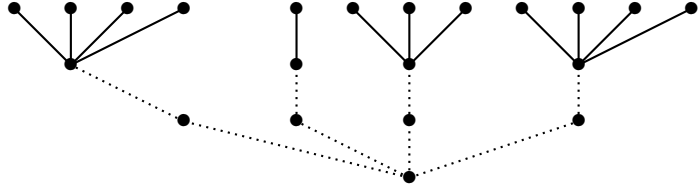

Next, we introduce the notion of a {\em shrubbery} of a tree poset $P$. 
Clearly, every leaf in $P$ belongs to a unique shrub in $P$. 
For example, consider the tree poset in Figure~\ref{F:Example}.
The tree poset in that figure has 4 subshrubs, whose Hasse diagrams are drawn in solid black lines. 
The {\em shrubbery of $P$} is the collection of subshrubs of $P$ that contain the leaves of $P$.

\begin{Theorem}\label{T:mainresult}
Let $P=\{\vep_1,\dots, \vep_m\}$ be a tree poset whose shrubbery consists of the shrubs, $S_{m_1},\dots, S_{m_s}$. 
Then the minimum distance of the code $\mc{C}_{\mb{O}_P}$ is given by 
\[
d(\mc{C}_{\mb{O}_P}) = (q-1)^{m- \sum_{i=1}^s (m_i-1)} (q-2)^{\sum_{i=1}^s (m_i-1)}.
\]
\end{Theorem}
\begin{proof}
By Proposition~\ref{P:pyramid}, the minimum distance $\mc{C}_{\mb{O}_P}$ is equal to $(q-1) d(\mc{C}_{\mb{O}_{P'}})$,
where $P'$ is the rooted forest obtained from $P$ by removing $\vep_1$. Let $P_1,\dots, P_r$ denote the connected components of $P'$. Then each $P_i$ ($i\in \{1,\dots, r\}$) is a rooted tree.
By repeatedly applying Proposition~\ref{P:dofproducts} and Proposition~\ref{P:pyramid}, we reach to the shrubberies of the $P_i$'s for all $i\in \{1,\dots, r\}$. The union of the shrubberies of the $P_i$'s ($i\in \{1,\dots, r\}$) is equal to the shrubbery of $P$, that is, 
$S_{m_1},\dots, S_{m_s}$.  
For $l\in \{1,\dots,s\}$, the index $m_l$ is the number of vertices in the shrub $S_{m_l}$.  
Let $j$ denote the difference $m- \sum_{l=1}^s m_l$, which is equal to the number of vertices that are removed from $P$ to reach to the shrubbery $S_{m_1},\dots, S_{m_s}$.
In particular, we have the following formula for the minimum distance,
\begin{align}\label{A:beforesubs}
d(\mc{C}_{\mb{O}_P} )= (q-1)^j d(\mc{C}_{\mb{O}_{S_{m_1}}}) \cdot \hdots \cdot d(\mc{C}_{\mb{O}_{S_{m_s}}}).
\end{align}
We now observe that, for each $l\in \{1,\dots, s\}$, the order polytope $\mb{O}_{S_{m_l}}$ is a pyramid over the unit cube of dimension $m_l-1$. 
Therefore, by~\cite[Corollary 3.4]{SS2010}, the minimum distance of the corresponding code is given by $(q-1)(q-2)^{m_l-1}$.
Thus, by substituting these into (\ref{A:beforesubs}) we obtain the asserted formula for the minimum distance $\mc{C}_{\mb{O}_P}$. 
\end{proof}

\section{A Lemma on Ordinal Sums}\label{S:Main2}

Let $P$ and $Q$ be two posets. 
The {\em ordinal sum of $P$ and $Q$}, denoted by $P\oplus Q$, is the poset defined on the disjoint union $P\sqcup Q$ as follows.
Let $a$ and $b$ be two elements from $P\sqcup Q$. 
Then 
\begin{align*}
a\leq b \iff
\begin{cases}
\text{ if both of $a$ and $b$ are the elements of $P$, and $a\leq b$ in $P$};\\
\text{ if both of $a$ and $b$ are the elements of $Q$, and $a\leq b$ in $Q$};\\
\text{ if $a\in P$ and $b\in Q$}.
\end{cases}
\end{align*}

The order polytope of the ordinal sum of two posets can be described in terms of the order polytope of the summands. 
This relationship is expressed by the action of the {\em group of affine transformations} of a lattice.
To explain, let $\Z^k$ be a lattice, let $u$ be an element of $\Z^k$, and let $M$ an element of $GL_k(\Z)$. 
The map $T_{M,u}:\Q^k\to \Q^k$, defined by the formula $T(v) := M\cdot v + u$ for $v\in \Z^k$, is called an {\em affine transformation of $\Z^k$}.
Now, two polytopes $\mb{P}$ and $\mb{Q}$ in $\Z^k\otimes_\Z \Q \cong \Q^k$ are called {\em lattice equivalent} if there exists an affine transformation $T_{M,u}: \Q^k\to \Q^k$ such that $T_{M,u}(\mb{P}) = \mb{Q}$. 
Since the affine transformations form a group, the lattice equivalence is an equivalence relation on the collection of all polytopes in $\Q^k$. 
An important fact regarding the lattice equivalence is that  
two toric codes that are obtained from two lattice equivalent polytopes have the same parameters. 
For a detailed explanation of this fact, we refer the reader to~\cite[Section 4]{LS2007}.

\begin{Lemma}\label{L:ordinalfree}
Let $P$ and $Q$ be two posets. 
Then the order polytope of the ordinal sum $P\oplus Q$ is lattice equivalent to the free sum of polytopes $\mb{O}_P \oplus \mb{O}_{Q}$. 
\end{Lemma}

\begin{proof}
Let $n$ and $m$ denote the cardinalities of $P$ and $Q$ respectively. 
Then $\mb{O}_P\subset \Q^n$ and $\mb{O}_Q \subset \Q^m$. 
Let $I$ (resp. $I'$) be an element of $J(P^{opp})$ (resp. of $J(Q^{opp})$). 
By abuse of notation, we will use the same notation $I$ (resp. $I'$) for the upper order ideal generated by $I$ (resp. $I'$) in $P\oplus Q$.
In this notation, clearly, for every upper order ideal $I$ of $P$ we have $Q\leq I$ in $J( (P\oplus Q)^{opp})$. 
In terms of cartesian coordinates on $\Q^n\times \Q^m$, this fact amounts to the fact that $\rho_{P\oplus Q}(I)$ has 1's on its last $m$ coordinates.
In other words, in $\Q^n\times \Q^m$, the vector $v_0:=(0,\dots, 0,1,\dots, 1)$ corresponds to both of 
1) the empty upper order ideal of $P$, 2) the maximal upper order ideal of $Q$.
We now consider the affine translate $\mb{O}_{P\oplus Q}-v_0$ in $\Q^n\times \Q^m$. 
Under this translation, the vertices that correspond to the upper order ideal in $P$ are mapped to the negatives of the lower order ideals in $P$. 
Therefore, we have the following equality of polytopes: 
\begin{align*}
\mb{O}_{P\oplus Q} - v_0 = (-\mb{O}_{P^{opp}})\oplus \mb{O}_{Q}.
\end{align*}
But the polytope $-\mb{O}_{P^{opp}}$ is lattice equivalent to $\mb{O}_P$, hence, we obtain the equivalence,
\begin{align*}
\mb{O}_{P\oplus Q} - v_0 \cong \mb{O}_{P}\oplus \mb{O}_{Q}.
\end{align*}
This finishes the proof of our assertion.
\end{proof}

Recall that the minimum distance of the toric code that is obtained from the direct product of two polytopes $\mb{P}$ (in $\Q^m$) and $\mb{Q}$ (in $\Q^n$)
is given by the product of the minimum distances of the codes that are associated with $\mb{P}$ and $\mb{Q}$ (Lemma~\ref{L:SS2010T2.1}).
Let $h$ be a polynomial from $H^0(T_N(\F_q), \mb{P})$. 
The {\em weight of $h$}, denoted $wt(h)$, is the maximum number of nonzero coordinates in the image vector of the evaluation of $h$ on the points of $T_N(\F_q)$.
Let $f$ be a polynomial from $H^0(T_N(\F_q), \mb{P})$ such that $wt(f) = d(\mc{C}_{\mb{P}})$. 
Similarly, let $g$ be a polynomial from $H^0(T_{N'}(\F_q), \mb{Q})$ such that $wt(g)=d(\mc{C}_{\mb{Q}})$.
In their proof of Lemma~\ref{L:SS2010T2.1}, Soprunov and Soprunova~\cite[Theorem 2.1]{SS2010} show that the weight of the polynomial $fg$ is equal to $d(\mc{C}_{\mb{P}\times \mb{Q}})$. 
Note that $f$ and $g$ separately belong also to the space of sections $H^0(T_{N\times N'}(\F_q) , \mb{P}\oplus \mb{Q})$. 
This in particular gives us an upper bound for $d(\mc{C}_{ \mb{P}\oplus \mb{Q}})$ as follows. 
Clearly, the total number of points in $T_{N\times N'}(\F_q)$ ($\cong (\F_q^*)^{m+n}$) where $f$ (resp. $g$) vanishes is given by $Z(f)(q-1)^n$ (resp. by $Z(g)(q-1)^m$). 
Thus, we have 
\begin{align*}
d(\mc{C}_{ \mb{P}\oplus \mb{Q}}) \leq \max \{ (q-1)^{m+n} - Z(f) (q-1)^n, (q-1)^{m+n} - Z(g) (q-1)^m\}. 
\end{align*}
Next, we apply this observation to an ordinal sum of posets.

Let $m$ be a positive integer. 
Let us denote an antichain with $m$ elements by $A_m$. 
The order polytope of $A_m$ is the $m$-dimensional unit cube. 
Note that an $m$-chain is given by $A_1\oplus \cdots \oplus A_1$ ($m$ copies), which we denote by $C_m$.

\begin{Lemma}\label{L:prepbipartite}
Let $m$ be a positive integer. 
Then the minimum distance of the toric code associated with $\mb{O}_{A_m\oplus A_m}$ is given by $(q-1)^{m}(q-2)^m$.
\end{Lemma}

\begin{proof} 
We begin with a slightly more general setup. 
Let $m \leq n$ be two positive integers. We consider the ordinal sum $A_m\oplus A_n$. 
In the light of Lemma~\ref{L:ordinalfree}, we may assume that $\mb{O}_{A_m\oplus A_n} = \mb{O}_{A_m}\oplus \mb{O}_{A_n}$. 
Let $f$ be a polynomial in $H^0(T_N(\F_q), \mb{O}_{A_m})$ such that $wt(f) = d(\mc{C}_{\mb{O}_{A_m}})$. 
Then we know that 
\[
Z(f) = (q-1)^m - d(\mc{C}_{\mb{O}_{A_m}})= (q-1)^m - (q-2)^m.
\]
Similarly, let $g$ be a polynomial in $H^0(T_{N'}(\F_q), \mb{O}_{A_n})$ such that $wt(g) = d(\mc{C}_{\mb{O}_{A_n}})$. 
Then we know that 
\[
Z(g) = (q-1)^n - d(\mc{C}_{\mb{O}_{A_n}})= (q-1)^n - (q-2)^n.
\]
Therefore, the minimum distance of $\mb{O}_{A_m\oplus A_n}$ is bounded by 
\begin{align*}
d(\mc{C}_{ \mb{O}_{A_m}\oplus \mb{O}_{A_n}})  \leq \max \{ &(q-1)^{m+n}  - ((q-1)^m - (q-2)^m)  (q-1)^n,\\ &(q-1)^{m+n} - ((q-1)^n - (q-2)^n)  (q-1)^m\} \\
&= \max \{ (q-2)^m (q-1)^n, (q-2)^n (q-1)^m \} \\ 
&= (q-2)^m(q-1)^n.
\end{align*}
In particular, if $m=n$, then we see that 
\begin{align}\label{A:firstbound}
d(\mc{C}_{ \mb{O}_{A_m}\oplus \mb{O}_{A_n}}) \leq (q-2)^m (q-1)^m.
\end{align}
We notice that the poset $A_m\oplus A_m$ is covered by $m$ 2-chains, $H_m:= \sqcup_{i=1}^m C_2$. 
It is easy to check the polytope containment 
\[
\mb{O}_{A_m\oplus A_m} \subseteq \mb{O}_{H_m}.
\]
This means that the space of sections of the line bundle determined by $\mb{O}_{A_m\oplus A_m}$ is contained in the space of sections of the line bundle determined by $\mb{O}_{H_m}$.
Since these sections are evaluated on the same torus, the minimum distance of the code $\mc{C}_{ \mb{O}_{A_m}\oplus \mb{O}_{A_m}}$ is bounded from below by the minimum distance of $\mc{C}_{\mb{O}_{H_m}}$, which is equal to $(q-1)^m (q-2)^m$. 
The rest of the proof follows from (\ref{A:firstbound}).
\end{proof}

\begin{Theorem}\label{T:mainresult2}
Let $m$ be a positive integer.
The minimum distance of a toric code associated with an $(m,m)$-bipartite poset is given by $(q-1)^m(q-2)^m$. 
\end{Theorem}

\begin{proof}
Let $H_m$ denote $\sqcup_{i=1}^m C_2$.
By the proof of Lemma~\ref{L:prepbipartite}, we know that 
\[
d(\mc{C}_{\mb{O}_{A_m\oplus A_m}}) = d(\mc{C}_{\mb{O}_{H_m}}) = (q-1)^m (q-2)^m.
\]
It is easy to check (by computing the vertices of the order polytopes) that if $P$ is an $(m,m)$-bipartite poset, then 
$\mb{O}_{A_m\oplus A_m} \subseteq \mb{O}_P \subseteq {\mb{O}_{H_m}}$. 
These inclusions give the following inequalities:
\[
d(\mc{C}_{\mb{O}_{A_m\oplus A_m}}) \geq d(\mc{C}_{\mb{O}_P})\geq  d(\mc{C}_{\mb{O}_{H_m}}),
\]
which are actually equalities. This finishes the proof of our theorem.
\end{proof}

\begin{Proposition}
Let $m$ be a positive integer. 
Then we have the following formulas for the dimensions of the toric codes associated with $A_m\oplus A_m$ and $H_m:=\sqcup_{i=1}^m C_2$. 
\begin{enumerate}
\item $\dim \mc{C}_{\mb{O}_{A_m\oplus A_m}} = 2^{m+1}-1$, and 
\item $\dim \mc{C}_{\mb{O}_{H_m}} = 3^m$.
\end{enumerate}
\end{Proposition}
\begin{proof}
The dimension of a toric code defined by an order polytope is equal to the number of vertices of the polytope.
In the former case, we have the free sum of two $m$ dimensional cubes.
Therefore, the dimension in this case is given by $2^m + 2^m -1 = 2^{m+1}-1$. 
In the latter case, the vertices of $\mb{O}_{H_m}$ are given by the upper order ideals in $H_m$. 
Any such idea is uniquely determined by $a$ minimal elements $\hat{0}_{i_1},\dots, \hat{0}_{i_a}$ in $H_m$, and $b$ 
maximal elements $\hat{1}_{j_1},\dots, \hat{1}_{j_b}$, where $\hat{1}_{j_r}$ ($1\leq r \leq b$) does not cover any element from $\{\hat{0}_{i_1},\dots, \hat{0}_{i_a}\}$.
Therefore, the total number of such upper order ideals is given by $\sum_{a= 0}^m \sum_{b=0}^{m-a} {m\choose a} {m-a\choose b}$. 
By using the binomial theorem, we see that this sum is equal to $3^m$. 
\end{proof}

\begin{Example} 

We consider the posets $P_1,P_2$ and $P_3$ that are defined in Figure~\ref{F:N}.
In Table~\ref{T:Jopp} we listed their upper order ideals. 

\begin{figure}[htp]
\centering
\begin{tikzpicture}

\begin{scope}[xshift = -3cm]
\node (a) at (-.5,0) {$\vep_1$};
\node (b) at (.5,0) {$\vep_2$};
\node (c) at (-.5,1) {$\vep_3$};
\node (d) at (.5,1) {$\vep_4$};
\draw[-, thick] (a) to (c);
\draw[-, thick] (a) to (d);
\draw[-, thick] (b) to (c);
\draw[-, thick] (b) to (d);
\end{scope}

\begin{scope}
\node (a) at (-.5,0) {$\vep_1$};
\node (b) at (.5,0) {$\vep_2$};
\node (c) at (-.5,1) {$\vep_3$};
\node (d) at (.5,1) {$\vep_4$};
\draw[-, thick] (a) to (c);
\draw[-, thick] (b) to (c);
\draw[-, thick] (b) to (d);
\end{scope}

\begin{scope}[xshift = 3cm]
\node (a) at (-.5,0) {$\vep_1$};
\node (b) at (.5,0) {$\vep_2$};
\node (c) at (-.5,1) {$\vep_3$};
\node (d) at (.5,1) {$\vep_4$};
\draw[-, thick] (a) to (c);
\draw[-, thick] (b) to (d);
\end{scope}

\end{tikzpicture}
\caption{The posets $P_1,P_2$, and $P_3$ (from left to right).}
\label{F:N}
\end{figure}
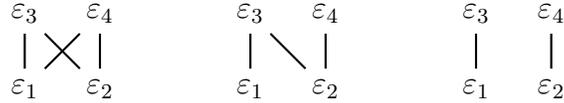

\begin{table}
\begin{center}
\begin{tabular}{ c|c|c } 
$J(P^{opp}_1)$ & $J(P^{opp}_2)$  & $J(P^{opp}_3)$  \\
\hline
$\{\vep_1,\vep_2,\vep_3,\vep_4\}$ & $\{\vep_1,\vep_2,\vep_3,\vep_4\}$ & $\{\vep_1,\vep_2,\vep_3,\vep_4\}$  \\
\hline
$\{\vep_1,\vep_3,\vep_4\}$ & $\{\vep_1,\vep_3,\vep_4\}$ & $\{\vep_1,\vep_3,\vep_4\}$  \\
\hline
$\{\vep_2,\vep_3,\vep_4\}$ & $\{\vep_2,\vep_3,\vep_4\}$ & $\{\vep_2,\vep_3,\vep_4\}$  \\
\hline
$\{ \vep_3,\vep_4\}$ & $\{ \vep_3,\vep_4\}$ & $\{\vep_3,\vep_4\}$  \\
\hline
$\{ \vep_3 \}$ & $\{ \vep_1,\vep_3\}$ & $\{\vep_1,\vep_3\}$  \\
\hline
$\{ \vep_4 \}$ & $\{ \vep_3\}$ & $\{\vep_2,\vep_4\}$  \\
\hline
$\emptyset$ & $\{ \vep_4\}$ & $\{\vep_3\}$  \\
\hline
 & $\emptyset$ & $\{\vep_4\}$  \\
\hline
 &  & $\emptyset$  \\
\hline
\end{tabular}
\end{center}
\caption{The upper order ideals of $P_1,P_2,P_3$.}
\label{T:Jopp}
\end{table}
The minimum distance of the toric code associated with the order polytope of $P_i$ ($i\in \{1,2,3\}$) equals 
\[
d(\mc{C}_{\mb{O}_{P_i}}) = (q-1)^2 (q-2)^2.
\]

\end{Example}

\section*{Acknowledgement}
We thank Roy Joshua and G.V. Ravindra for many useful discussions on the topics of this paper.
We thank the referee for the constructive comments and for the very careful reading of our paper. 
Finally, we gratefully acknowledge the research of the first author was partially supported by a grant from the Louisiana Board of Regents.

\bibliographystyle{plain}
\bibliography{references}

\end{document}